\documentclass[11pt]{article}
\usepackage{tipa}
\usepackage{stmaryrd}

\usepackage{color,amsmath,amssymb, amsfonts,amstext,amsthm}

\usepackage{epsfig, graphicx, graphics}
\usepackage{bbm}
\usepackage[colorlinks]{hyperref}
\usepackage{mathrsfs}
\usepackage{url}

\newcommand{\R}{\mathbb{R}}

\newcommand{\Rn}{\mathbb{R}^n}

\renewcommand{\mid}{:}

\newtheorem{theorem}{Theorem}

\newtheorem{definition}[theorem]{Definition}
\theoremstyle{remark}
\newtheorem{remark}{Remark}
\newtheorem{example}[remark]{Example}
\topsep=\parskip

\title{Small noise approximation of center manifolds for stochastic dynamical systems
\footnote{This work was partly supported by the NSF Grant  1025422,   and the NSFC grants 10971225 and 11028102. } }

\author{Jian Ren, Zhongkai Guo, Xianming Liu and Xiangjun Wang
\\School of Mathematics and Statistics\\
Huazhong University of Science and Technology, Wuhan, 430074, China \\
 \& \\Institute for Pure and Applied Mathematics, University of California, Los Angeles, CA 90095, USA
 \\Email:   renjian0371@gmail.com   }

\date{\today}

\begin{document}
\maketitle

\begin{abstract}
 This paper provides a small noise approximation for local random center manifolds
 of a class of stochastic dynamical systems in Euclidean space.
 An example is presented to illustrate the method.

\end{abstract}




\bigskip
\section{Introduction}

Center manifolds, together with stable and unstable manifolds, provide geometric structures that help understand nonlinear and stochastic dynamics. For stochastic dynamical systems, it is difficult to visualize or depict center manifolds. In the paper, we derive a small noise approximation for random center manifolds.

We consider a system of nonlinear random differential equations
(RDE)   on a probability space $(\Omega,\mathcal{F}, \mathbb{P} )$,
\begin{equation}\label{eq(2.1)}
    \frac{d u}{dt}=A(\theta_t \omega) u + F(\theta_t \omega, u),
\end{equation}
where $u\in \Rn$,   $\theta_t: \Omega \to \Omega$ is an ergodic flow
that preserves the probability measure $\mathbb{P}$, $A$ is a
$n\times n$ matrix satisfies integrability condition, and $F$ is a
nonlinear mapping in $\Rn$. Assume that $0$ is a fixed point, i.e.,
$F(\theta_t \omega, 0)=0$ for all $t$ and all $\omega$. We also
assume that the Jacoboian matrix $D_u F(\theta_t \omega,
0)=\textbf{0}$, which implies that $F$ is indeed nonlinear.

This RDE system is often coming from a stochastic differential
equation (SDE) system, by a stationary coordinate
transform \cite{Duan}.

Boxler \cite{Boxler, Boxler2} proved a random center manifold
theorem, and also considered reduction to the center manifold. Power
series, in terms of variables on the center manifold, approximations
are also available  in \cite{Boxler2, Xu}.  We also mention some earlier explorations on center manifolds for stochastic systems  \cite{Knobloch, Schumaker}.

 In this present paper,
 we derive a small noise approximation of random center manifolds, for the stochastic system (\ref{eq(2.1)}) where the
 linear operator $A(\theta_t \omega)$ has an exponential trichotomy property.

\section{Exponential trichotomy}

Let $(\Omega,\mathcal{F}, \mathbb{P})$ be a probability space  with an ergodic flow $\theta_t$.
A measurable map $\phi(t, \omega, u):\mathbb{R}\times  \Omega\times
\mathbb{R}^n \rightarrow \mathbb{R}^n$   is called a random dynamical system (RDS) if it
satisfies the following conditions
\begin{enumerate}
\item   $\phi(0,\omega, u)= u$;
\item    $\phi(t+s, \omega, u)=\phi(t, \theta_s\omega, u)\phi(s, \omega, u)$,   for all $t, s \in \mathbb{R}$ and    all $\omega \in \Omega$.
\end{enumerate}
The latter  is the so-called cocycle property.

\bigskip


For the linear RDE system in $\Rn$
\begin{equation}\label{linear}
    \frac{d u}{dt}=A(\theta_t \omega) u ,
\end{equation}
 where the random matrix $A(\theta_t \omega)$ satisfies an integrability
 condition, i.e.
 $A\in \mathbf{L}^1(\Omega,\mathcal{F}, \mathbb{P})$. The
integrability of $A(\theta_t\omega) $ implies that the RDE system
\eqref{linear} generates a linear cocycle $\Phi(t, \omega)$, which
is the fundamental matrix  $\Phi(t, \omega) \triangleq
e^{\int_0^tA(\theta_\tau\omega)\,\mathrm{d}\tau}$.
Furthermore, this linear cocycle $\Phi(t, \omega)$ satisfies an integrability condition for  the multiplicative ergodic theorem
(MET) in \cite[Chapter 3]{Arn98}. This MET   yields that there exist Lyapunov
exponents $\lambda_1>\lambda_2
>\cdots>\lambda_r$, with $r\leq n$, together with $r$ Osledets
subspaces $E_i(\omega)$, so that $\Rn$ decomposes into the direct
sum:
\begin{eqnarray}\label{subspace decompose}
\Rn = E_1(\omega) \oplus E_2(\omega) \oplus \hdots E_r(\omega).
\end{eqnarray}
All the subspaces $E_i(\omega)$ are measurable and are random invariant
subspaces of $\Phi(t, \omega)$, i.e. for $i=1,2,\cdots, r$,
\begin{eqnarray}\label{invariant}
\Phi(t, \omega)E_i(\omega)=E_i(\theta_t\omega),\quad \text{for all}
\,\,\,t\in \mathbb{R}, \,\,\,\mathbb{P}-a.s.
\end{eqnarray}
We define stable, center and
unstable subspaces
\begin{eqnarray}
E^s(\omega) \triangleq \oplus_{\lambda_i<0} E_i(\omega), \\
E^c(\omega) \triangleq \oplus_{\lambda_i = 0} E_i(\omega), \\
E^u(\omega) \triangleq \oplus_{\lambda_i>0} E_i(\omega).
\end{eqnarray}
Therefore,  $$\Rn = E^s(\omega) \oplus E^c(\omega) \oplus
E^u(\omega),$$ and $\Phi(t, \omega)$ can be decomposed as $$\Phi(t,
\omega)=\Phi_s(t, \omega)\oplus\Phi_c(t, \omega)\oplus\Phi_u(t,
\omega).$$
 Denote
$$
\lambda_s \triangleq \max_{\lambda_i <0} \lambda_i<0 ,\quad
\lambda_u \triangleq \min_{\lambda_i >0} \lambda_i>0.
$$
Then
\begin{eqnarray*}
\lim\limits_{t\rightarrow\pm\infty}\frac{1}{t}\log|\Phi(t,
\omega)v|&=&\lambda ,
\end{eqnarray*}
with
$$
\lambda = \begin{cases}
 \lambda_s,\;\; \text{for}\;\; v\in E^s(\omega)\setminus\{0\},\\
 0,\;\;\;\; \text{for} \;\;    v\in E^c(\omega)\setminus\{0\}, \\
 \lambda_u,\;\; \text{for}\;\; v\in E^u(\omega)\setminus\{0\}.
         \end{cases}
$$
 Taking a positive number $\gamma < \frac12\min\{|\lambda_s|, \lambda_u \} $,
    the MET leads to the following estimation (see \cite[Lemma 4.1]{Boxler} or
    \cite[p.143]{Boxler2}):
  \begin{eqnarray}\label{exponential estimate}
|\Phi_s(t, \omega)v|&\leq& K_s(\omega) e^{(\lambda_s+\gamma) t}| v |, \quad v\in
E^s,\quad t\geq 0,\nonumber\\
|\Phi_c(t, \omega)v|&\leq& K_c(\omega)  e^{\gamma |t|}| v |, \quad \;\;\quad v\in
E^c,\quad t\in
\mathbb{R},\nonumber\\
|\Phi_u(t, \omega)v|&\leq&  K_u(\omega)  e^{(\lambda_u-\gamma) t}| v |, \quad v\in
E^u,\quad t \leq 0,
\end{eqnarray}
where $K_s, K_c, K_u$ are   measurable functions depending on $\gamma$ (with values in $[1, \infty)$).

\bigskip
 Similar to the definition of exponential dichotomy in
\cite{Duan}, the cocycle $\Phi(t, \omega)$  with the conditions
\eqref{subspace decompose}, \eqref{invariant} and \eqref{exponential
estimate}   is said to have   exponential trichotomy.

\bigskip

For the nonlinear RDE system
\begin{equation}\label{nonlinear}
    \frac{d u}{dt}=A(\theta_t \omega) u + F(\theta_t\omega, u),
\end{equation}
we assume that the nonlinear term~$F$  is
  Lipschitz continuous on~$\mathbb{R}^n$, that is,
\[
|F(u_1)-F(u_2)|\le \L_F\,|u_1-u_2 |,
\]
with the      Lipschitz constant $\L_F>0$. In the next section, we
consider the random center manifold for this system.

If  the  nonlinear term~$F$ is only locally Lipschitz, then we let
$F^{(R)}(u)=\chi_R(u)F(u)$, where $\chi_R(u)$~is a cut-off function.
Thus  $F^{(R)}$~is  globally Lipschitz with Lipschitz
constant~$R\L_F$; see \cite[Lemma 4.1]{Carr}. In this case we would
obtain a local random center manifold.

By using the fundamental function
\begin{equation*}
\Phi(t, \omega)=e^{\int_0^tA(\theta_\tau\omega)\,\mathrm{d}\tau},
\end{equation*}
 the solution of  the random differential  system   \eqref{nonlinear}   can be interpreted as
\begin{equation*}\label{sol-nonlinear}
u(t, \omega, u_0)=\Phi(t, \omega)u_0 +\int_0^t \Phi(t-\tau,
\omega)F(\theta_\tau\omega, u(\tau))\,\mathrm{d}\tau\,.
\end{equation*}
This equation has  a unique measurable solution due to the Lipschitz
continuity of $F$, and the solution mapping $(t, \omega, u_0)\mapsto
u(t, \omega, u_0)$ generates a random dynamical system $\phi(t,
\omega, u_0)$. Moreover, the Jacobian mapping $D_{u_0}u(t, \omega,
u_0)=~\Phi(t, \omega)$ satisfies the exponential trichotomy.





\section{Random center manifold}  \label{small888}

In this section, we consider random center manifolds for a SDE
system with small multiplicative     noise. We will recall a random
center manifold theorem  \cite{Boxler} and then approximate the
random center manifold under small noise, in the case of
non-positive Lyapunov exponents. This is done in deterministic case
in \cite{Carr} and \cite[Chapter 3]{GH}.

\bigskip

We introduce the definition of random center manifolds.
\begin{definition}[Random Center Manifold]
A random set~$M(\omega)$ is called a  random center manifold for a random dynamical system ~$\phi(t, \omega, x)$, if
 it satisfies the following conditions

(i) It is an   invariant set, i.e.  $\phi(t, \omega,
M(\omega))\subset~ M(\theta_t\omega)$ for all $t$.

(ii) It can be represented as a graph of a
  mapping from the center subspace to its complement, i.e.
there is a mapping $h^c(\omega, \cdot): E^c\to E^s\oplus E^u$, such
that $M(\omega)=\big\{\big(v, h^c(\omega, v)\big) \mid v\in
E^c\big\}$, where $h^c(\omega, 0)=0$,  $h^c(\cdot, v)$ is measurable
for every $v\in E^c$ and the tangency condition
$Dh^c(\omega, 0)=0$ holds. The center manifold $M(\omega)$ is often denoted as~$M^c(\omega)$.\\
\\
It is called a Lipschitz center manifold if the  mapping $h^c(\omega, \cdot): E^c\to E^s\oplus E^u$ is Lipschitz and   the tangency condition is absent.
We call $M^c(\omega)$  a local center manifold if it is a graph of the
   mapping~$\chi_R(v)\;h^c(\omega, \cdot)$, where $\chi_R(v)$ is a cut-off function.
\end{definition}

\bigskip

For the nonlinear SDE system defined in the previous section,
  Boxler  \cite{Boxler, Boxler2} has shown that there exists a
   (local) random center manifold, as the graph of a mapping $h^c$ that satisfies the following Liapunov-Perron integral equation
   \begin{eqnarray} h^c(\omega,
x^c)&=&\bigg(\int_{-\infty}^0\Phi^{-1}_s(\tau,
\omega)F^s\Big(\theta_\tau\omega, \big( \phi^c(\tau, \omega, x^c, h^c(\omega, x^c)), h^c(\theta_\tau\omega, \phi^c(\tau, \omega, x^c, h^c(\omega, x^c))) \big) \Big)\,\mathrm{d}\tau,\nonumber\\
&&{}\int_{\infty}^0\Phi^{-1}_u(\tau,
\omega)F^u\Big(\theta_\tau\omega, \big(\phi^c(\tau, \omega, x^c,
h^c(\omega, x^c)), h^c(\theta_\tau\omega, \phi^c(\tau, \omega, x^c,
h^c(\omega, x^c)))
\big) \Big)\,\mathrm{d}\tau\bigg), \nonumber\\
\end{eqnarray}
for $x^c\in E^c$ and $\omega \in \Omega$ where $F^s$, $F^u$ are
  the projection of $F$ to the stable and unstable
   subspaces $E^s, E^c$ and $E^u$, respectively, and $\phi^c$ is the projection of $\phi$ to the center
  subspace $E^c$.


\subsection{Center manifold of  a system with small  noise}
\label{center-small para}

 We consider a SDE system in Stratonovich form
\begin{eqnarray}
\begin{cases}\label{stochastic center}
\dot x = A^c x + f^c(x, y) + \big(\varepsilon  x^{\intercal}\circ\dot W_t^1\big)^\intercal,\,\,\, x\in R^n, \\
\dot y = A^s y + f^s(x, y) + \big(\varepsilon y^{\intercal}\circ\dot
W_t^2\big)^\intercal,\,\,\,y\in R^m,
\end{cases}
\end{eqnarray}
where $^\intercal$ indicates the transpose of a vector or a matrix.
In this system,   $A^c$ and  $A^s$ are respectively $n\times n$ and
$m\times m$ matrixes. The   nonlinear functions $f^c(x,
y):\mathbb{R}^n\times\mathbb{R}^m\rightarrow\mathbb{R}^n$ and
$f^s(x, y):\mathbb{R}^n\times\mathbb{R}^m\rightarrow\mathbb{R}^m$
are $C^1$-smooth and Lipschitz with Lipschitz constant $L_f$, i.e.
$f^c(x, y)$ and $f^s(x, y)$ satisfy
$$|f^c(x_1, y_1) - f^c(x_2,
y_2)|\leq L_f\,\big(|x_1 -x_2|+|y_1-y_2|\big), $$and
$$|f^s(x_1, y_1) -f^s(x_2, y_2)|\leq L_f\,\big(|x_1
-x_2|+|y_1-y_2|\big),
$$
where $|x|$ is the Euclidean norm of $x\in\mathbb{R}^n $ or
$x\in\mathbb{R}^m$. Usually, $f^c, f^s$ are locally Lipschitz
continuous and in that case we would get a local random center
manifold.  The noise intensity $\varepsilon$ is a small positive
parameter. Moreover, $\{W_t^i:t\in\mathbb{R},\,\,i=1,2\}$ are
$n\times n$ and $m\times m$ matrixes, respectively,  with two-sided
scalar Wiener process $W_t$ as principal diagonal elements and all
other elements being zero.

 For matrixes $A^c$ and $A^s$, we make the
following assumption.

$\mathbf{(H)}$: There are positive constants $\beta$, $\gamma$ and
$K$, satisfying $\beta>\gamma\geq0$ and $K>0$,  such that for every
$x\in R^n$ and $y\in R^m$, the following exponential estimates
hold
 $$|e^{A^ct}x|_{R^n} \leq K e^{\gamma|t|}|x|_{R^n},\quad t\in \mathbb{R};\quad \quad |e^{A^st}y|_{R^m} \leq
 K e^{-\beta t}|y|_{R^m},\quad t \geq 0.$$

To facilitate random dynamical systems approach, we convert
\eqref{stochastic center} into a system of random differential
equations (RDE). To this end, consider   linear stochastic differential
equations
\begin{eqnarray}\label{O-U}
 \mathrm{d}z_i + z_i \mathrm{d}t = \mathrm{d}W_t^i,\quad i=1,2.
\end{eqnarray}
Each of these equations has a   unique stationary solution
$$z_i(\omega)=\int_{-\infty}^0 e^{\tau}\,\mathrm{d}W_\tau^i,\quad i=1,2.$$
Moreover,
$$z_i(\theta_t\omega)=\int_{-\infty}^t
e^{\tau-t}\,\mathrm{d}W_\tau^i,\quad i=1,2.$$ Thus, the linear
differential equations,
\begin{eqnarray}\label{O-U-eps}
 \mathrm{d}Z_i + Z_i \mathrm{d}t = \varepsilon \mathrm{d}W_t^i,\quad i=1,2.
\end{eqnarray}
have unique stationary solutions $Z_i(\omega)=\varepsilon
z_i(\omega)$ and $Z_i(\theta_t\omega)=\varepsilon
z_i(\theta_t\omega)$ for $i=1,2$.

 Note that $z_i(\omega)$, $z_i(\theta_t\omega)$, $Z_i(\omega)$ and $Z_i(\theta_t\omega)$   for
$i=1,2$ are all principals diagonals matrixes, with the principal
diagonals values $\int_{-\infty}^0 e^{\tau}\,\mathrm{d}W_\tau$,
$\int_{-\infty}^t e^{\tau-t}\,\mathrm{d}W_\tau$,
$\varepsilon\int_{-\infty}^0 e^{\tau}\,\mathrm{d}W_\tau$ and
$\varepsilon\int_{-\infty}^t e^{\tau-t}\,\mathrm{d}W_\tau$
respectively.

\bigskip

 We introduce a random transformation
\[ \begin{matrix}\begin{pmatrix}  X \\ Y \\\end{pmatrix}
:=\mathcal V_\varepsilon(\omega,x, y)=
\begin{pmatrix}   e^{-\varepsilon z_1(\omega)}x \\ e^{-\varepsilon
z_2(\omega)}y
\\ \end{pmatrix} \end{matrix}.\]
Under this transformation, the SDE system \eqref{stochastic center}
is converted into a RDE system
\begin{eqnarray}
\begin{cases}\label{random center}
\dot X = A^c X + \varepsilon z_1 X + F^c(X, Y),\,\,\,X\in R^n,\\
\dot Y = A^s Y + \varepsilon z_2 Y + F^s(X, Y),\,\,\,Y\in R^m.
\end{cases}
\end{eqnarray}
with
$$F^c(X, Y) = e^{-\varepsilon z_1(\theta_t\omega)}
f^c(e^{\varepsilon z_1(\theta_t\omega)} X,  e^{\varepsilon
z_2(\theta_t\omega)} Y),
$$
\\and
$$F^s(X, Y) = e^{-\varepsilon z_2(\theta_t\omega)} f^s(
e^{\varepsilon z_1(\theta_t\omega)}X, e^{\varepsilon
z_2(\theta_t\omega)} Y),
$$
being Lipschitz functions with the Lipschitz
constant $L_f$.

For every $\eta > 0$, 
 we define a product Banach Space $C_\eta:=C_\eta^1\times
 C_\eta^2$, where\\
$C_\eta^1 = \big\{\varphi:\mathbb{R}\rightarrow \mathbb{R}^{n}
 :\phi$  is continuous and $\sup\limits_{t\in \mathbb{R}} |\,e^{-\eta |t|}e^{-\varepsilon\int_0^t z_1(\theta_r\omega)\,\mathrm{d}r}\varphi(t)|_{\mathbb{R}^{n}}<
\infty\big\},$\\$C_\eta^2 = \big\{\varphi:\mathbb{R}\rightarrow
\mathbb{R}^{m}
 :\phi$  is continuous and $\sup\limits_{t\in \mathbb{R}} |\,e^{\eta t}e^{-\varepsilon\int_0^t z_2(\theta_r\omega)\,\mathrm{d}r}\varphi(t)|_{\mathbb{R}^{m}}<
\infty\big\},$\\with  norm
\[|\varphi(t)|_{C_\eta^1} = \sup_{t\in \mathbb{R}}|\,e^{-\eta |t|}e^{-\varepsilon\int_0^t
z_1(\theta_r\omega)\,\mathrm{d}r}\varphi(t)|_{\mathbb{R}^{n}},
 \] and \[|\phi(t)|_{C_\eta^2} = \sup_{t\in \mathbb{R}}|\,e^{\eta t}e^{-\varepsilon\int_0^t z_2(\theta_r\omega)\,\mathrm{d}r}\varphi(t)|_{\mathbb{R}^{m}}.
 \]
respectively. Furthermore, the norm in $C_\eta$  is
 $$ |(X, Y)|_{C_\eta} = |X|_{C_\eta^1} + |Y|_{C_\eta^2}.$$

The integral form of the differential equation \eqref{random
 center} can be written as
 \begin{eqnarray}
\begin{cases}\label{random solution}
X(t) = e^{A^ct + \varepsilon\int_0^t
z_1(\theta_r\omega)\,\mathrm{d}r} X(0) + \int_0^t
e^{A^c(t-\tau)+\varepsilon\int_\tau^t z_1(\theta_r\omega)\,\mathrm{d}r} F^c(X, Y)\,\mathrm{d}\tau,\\
Y(t) =\int_{-\infty}^t e^{A^s(t-\tau)+\varepsilon\int_\tau^t
z_2(\theta_r\omega)\,\mathrm{d}r}F^s(X, Y)\,\mathrm{d}\tau.
\end{cases}
\end{eqnarray}
Due to the Lipschitz condition of $F^c$ and $F^s$ in \eqref{random
center} and the measurability of $z_i(\omega)$, $i=1,2$. the RDE
system \eqref{random center} has a unique solution. Therefore, the
solution mapping generates a RDS. By the assumption (\textbf{H}) on
$A^c$, $A^s$ and the propersties of $z_i(\omega)$, $i=1,2$, as in
\cite{Duan}, the exponential trichotomy holds (with $E^u$ empty).

Boxler\cite{Boxler, Boxler2} obtained the existence of center
manifold by introducing random form. Here we have following result
by using Lyapunov-Perron method similar as in \cite{Duan}.

\begin{theorem}[]  Assume that $\mathbf{H}$ holds and that $\gamma<\eta<\beta$ satisfy the gap condition $\frac{KL_f}{\eta-\gamma}+\frac{KL_f}{\beta-\eta}<1$. Then there
exists a random center manifold $\mathcal{\tilde
H}^\varepsilon(\omega)=(\tilde\xi, \tilde H^\varepsilon(\omega,
\tilde \xi))$ for \eqref{random center}, where
$\tilde{H}^\varepsilon(\omega,
\cdot):~\mathbb{R}^n\rightarrow~\mathbb{R}^m$ is a Lipschitz
function.
\end{theorem}

\subsection{Small noise  approximation of random center manifolds} \label{approximation}

We now approximate the above random center manifold under small noise intensity $\varepsilon$. Denote the corresponding deterministic center
manifold (when $\varepsilon=0$) and the random center manifold  as
\begin{eqnarray}
\mathcal{\tilde H}^d(\tilde \xi)=\big\{\big(\tilde\xi, \tilde
H^d(\tilde\xi)\big)\big \},
\end{eqnarray}
and
\begin{eqnarray}
\mathcal{\tilde H}^\varepsilon(\omega, \tilde \xi
)=\big\{\big(\tilde\xi, \tilde H^\varepsilon(\omega, \tilde\xi
)\big)= \big(\tilde\xi, \tilde H^d(\tilde\xi) + \varepsilon \tilde
H^1(\omega, \tilde\xi) + \varepsilon^2 \tilde H^2(\omega,
\tilde\xi)+ \mathcal{O}(\varepsilon^3)\big) \big\},
\end{eqnarray}
 where
$\tilde H^d(\cdot)$, $\tilde H^1(\omega, \cdot)$ and $\tilde
H^2(\omega, \cdot)$ are locally defined on $\R^n$. Denoting
$x(0)=\xi\in\mathbb{R}^n$, as in \cite{Duan}, the center manifold
$\mathcal{H}^\varepsilon(\omega, \xi)$ can be converted back to the
center manifold for the original stochastic system \eqref{stochastic
center},
\begin{eqnarray}
\nonumber\mathcal{H}^\varepsilon(\omega, \xi) &=&\big\{\big(\xi,
H^\varepsilon(\omega, \xi)\big) : \xi\in\mathbb{R}^n\big\}= \mathcal
V_\varepsilon^{-1}\big( \cdot , \mathcal{\tilde
M}^\varepsilon(\omega)\big)\\&=&\big\{\big(e^{\varepsilon
z_1(\omega)}\tilde\xi, e^{\varepsilon z_2(\omega)}\tilde
H^\varepsilon(\omega, \tilde\xi)\big)\big\}= \big\{\big(\xi,
e^{\varepsilon z_2(\omega)}\tilde H^\varepsilon(\omega,
e^{-\varepsilon z_1(\omega)}\xi)\big)\big\}\nonumber\\&=&
\Big\{\Big(\xi, \tilde H^d(\xi)+\varepsilon \big(\tilde H^1(\omega,
\xi)+ z_2(\omega)\tilde H^d(\xi)- \tilde
H_\xi^d(\xi)z_1(\omega)\xi\big) + \nonumber\\&&{}\varepsilon^2
\big(\tilde H^d_\xi(\xi) z_1^2 \xi/2 -z_2 \tilde H^d_\xi(\xi) z_1
\xi + z_2^2\tilde H^d(\xi)/2 -\tilde H^1_\xi(\omega, \xi) z_1 \xi +
z_2\tilde H^1(\omega, \xi) + \tilde H^2(\omega, \xi)
\big)+\mathcal{O}(\varepsilon^3)\Big)\Big\}.\nonumber\\\label{inversion}
\end{eqnarray}
To approximate the center manifold, we expand
$$X(t)=X_0(t)+\varepsilon X_1(t)+\varepsilon^2 X_2(t)+~\cdots,\quad
\text{with}\,\,\, X(0)=\tilde\xi=X_0(0), \,\,\,
\tilde\xi\in\mathbb{R}^n,$$ and
\begin{equation}
Y(t)=Y_0(t)+\varepsilon Y_1(t)+\varepsilon^2 Y_2(t)+\cdots, \label{H=H^d+H^1}
\end{equation}
with
$$Y(0)=~\tilde H^\varepsilon(\omega, \tilde\xi)=~\tilde
H^d(\tilde\xi)+~\varepsilon \tilde H^1(\omega, \tilde\xi )+
\varepsilon^2 \tilde H^2(\omega,
\tilde\xi)+~\mathcal{O}(\varepsilon^3).$$

\bigskip
\noindent Noting that $$e^{-\varepsilon z_i} = 1 - \varepsilon z_i +
\frac{\varepsilon^2z_i^2}{2!} + \cdots,\quad \text{for}\,\,\,
i=1,2,$$ we have,
\begin{eqnarray*}
F^c\big(X(t), Y(t)\big) &=& e^{-\varepsilon z_1(\theta_t\omega)}
f^c\big(e^{\varepsilon z_1(\theta_t\omega)}X(t),  e^{\varepsilon
z_2(\theta_t\omega)}Y(t)\big) \\&=& f^c\big(X_0(t), Y_0(t)\big) +
\varepsilon\Big\{f^c_x\big(X_0(t), Y_0(t)\big)\big[X_1(t) +
z_1(\theta_t\omega) X_0(t)\big] \\
{}&&+ f^c_y\big(X_0(t), Y_0(t)\big)\big[Y_1(t) + z_2(\theta_t\omega)
Y_0(t)\big] - z_1(\theta_t\omega)f^c\big(X_0(t), Y_0(t)\big)\Big\}
\\{}&&+\varepsilon^2\Big\{z_1^2 f^c/2 + f_x^c\big[X_2-z_1^2 X_0/2  \big]
+ f_y^c\big[(z_2^2Y_0/2+z_2Y_1+Y_2)-
z_1(Y_1+z_2Y_0)\big]\\{}&&+f_{xx}^c(X_1+z_1X_0)^2/2+f_{xy}^c(X_1+z_1X_0)(Y_1+z_2Y_0)
+f_{yy}^c(Y_1+z_2Y_0)^2/2\Big\}+ \mathcal{O}(\varepsilon^3),
\end{eqnarray*}
and
\begin{eqnarray*}
F^s\big(X(t), Y(t)\big) &=& e^{-\varepsilon
z_2(\theta_t\omega)}f^s\big( e^{\varepsilon z_1(\theta_t\omega)}
X(t), e^{\varepsilon z_2(\theta_t\omega)} Y(t)\big)\\ &=&
f^s\big(X_0(t), Y_0(t)\big) + \varepsilon\Big\{f^s_x\big(X_0(t),
Y_0(t)\big)\big[X_1(t) +
z_1(\theta_t\omega) X_0(t)\big] \\
{}&& + f^s_y\big(X_0(t), Y_0(t)\big)\big[Y_1(t) +
z_2(\theta_t\omega) Y_0(t)\big] - z_2(\theta_t\omega)f^s\big(X_0(t),
Y_0(t)\big)\Big\} + \\{}&&\varepsilon^2\Big\{z_2^2f^s/2 +
f_x^s\big[X_2+z_1X_1+z_1^2X_0/2-z_2(X_1+z_1X_0)\big]+f_y^s(Y_2-z_2^2Y_0/2)\\{}&&+
f_{xx}^s(X_1+z_1X_0)^2/2 +f_{xy}^s(X_1+z_1X_0)(Y_1+z_2Y_0)+
f_{yy}^s(Y_1+z_2Y_0)^2/2\Big\}+\mathcal{O}(\varepsilon^3).
\end{eqnarray*}
Inserting these expansions into \eqref{random center} and matching
the terms with the same   powers of $\varepsilon$,  we get at the
$0$-th order

\[ \begin{matrix}\begin{pmatrix}  \dot X_0(t) \\ \dot Y_0(t) \\\end{pmatrix}
=\begin{pmatrix}  A^c \quad 0 \\ 0\quad A^s \\\end{pmatrix} \begin{pmatrix}X_0(t) \\  Y_0(t)\end{pmatrix}+\begin{pmatrix} f^c\big(X_0(t), Y_0(t)\big) \\
 f^s\big(X_0(t), Y_0(t)\big) \\\end{pmatrix}\end{matrix},\]
which can be expressed   as
\begin{eqnarray}\label{X0Y0}
\begin{pmatrix}   X_0(t) \\  Y_0(t) \\\end{pmatrix}
=\exp\bigg\{\begin{pmatrix}  A^c \quad 0 \\ 0\quad A^s
\\\end{pmatrix} t\bigg\}\begin{pmatrix}\tilde\xi \\\tilde
H_d(\tilde\xi)\end{pmatrix}+
\int_0^t \exp\bigg\{\begin{pmatrix}A^c\quad 0\\0 \quad A^s\end{pmatrix}(t-\tau)\bigg\}\begin{pmatrix} f^c\big(X_0(s), Y_0(s)\big) \\
f^s\big(X_0(t), Y_0(t)\big) \\\end{pmatrix}\,\mathrm{d}\tau.\nonumber\\
\end{eqnarray}
Furthermore, at the first order in $ \varepsilon$,
\begin{eqnarray}
\begin{matrix}\begin{pmatrix}  \dot X_1(t) \\ \dot Y_1(t)
\end{pmatrix}
&=&\bigg[\begin{pmatrix}  A^c \quad 0 \\ 0\quad A^s \\ \end{pmatrix} + \begin{pmatrix} f^c_x\big(X_0(t), Y_0(t)\big) \quad f^c_y\big(X_0(t), Y_0(t)\big) \\
 f^s_x\big(X_0(t), Y_0(t)\big)\quad f^s_y\big(X_0(t), Y_0(t)\big) \end{pmatrix}\bigg] \begin{pmatrix} X_1 \\
 Y_1 \\\end{pmatrix}\end{matrix}\nonumber \\{} + \begin{matrix}\begin{pmatrix}  f^c_x z_1(\theta_t\omega)X_0(t)+ z_1(\theta_t\omega)X_0(t) +
 f^c_yz_2(\theta_t\omega)Y_0(t) - z_1(\theta_t\omega)f^c \\f^s_xz_1(\theta_t\omega)X_0(t) + z_2(\theta_t\omega)Y_0(t) + f^s_yz_2(\theta_t\omega)Y_0(t) - z_2(\theta_t\omega)f^s \\\end{pmatrix}
 \end{matrix},\nonumber\\
\end{eqnarray}
which can be rewritten as
\begin{eqnarray}\label{X1Y1}
\nonumber\lefteqn{\begin{pmatrix}X_1(t)\\Y_1(t)\\\end{pmatrix}=
\exp\bigg\{\begin{pmatrix}A^c\quad 0\\0 \quad A^s\end{pmatrix}t+
\int_0^t\begin{pmatrix} f^c_x\big(X_0(s), Y_0(s)\big) \quad f^c_y\big(X_0(s), Y_0(s)\big) \\
 f^s_x\big(X_0(s), Y_0(s)\big)\quad f^s_y\big(X_0(s), Y_0(s)\big) \\\end{pmatrix}\,ds\bigg\}\begin{pmatrix}   X_1(0) \\  Y_1(0)
 \\\end{pmatrix}} \\&&+\int_0^t \exp\bigg\{\begin{pmatrix}A^c\quad 0\\0 \quad A^s\end{pmatrix}(t-\tau)-\int_s^t\begin{pmatrix} f^c_x \quad f^c_y \\
 f^s_x\quad f^s_y \\\end{pmatrix}\,dr\bigg\}\begin{pmatrix}  f^c_xz_1X_0 + z_1X_0 + f^c_yz_2Y_0 - z_1f^c \\f^s_xz_1X_0 + z_2Y_0 + f^s_yz_2Y_0 - z_2f^s \\\end{pmatrix}
\mathrm{d}\tau,\nonumber\\
 \end{eqnarray}
and when $\varepsilon$ is order 2,
\begin{eqnarray}\label{X2}
\dot X_2(t)&=&A^cX_2(t)+z_1X_1(t)+z_1^2f^c/2  +
f_y^c\big[z_2^2Y_0/2+z_2Y_1+Y_2- z_1(Y_1+z_2Y_0)\big]+
f_x^c(X_2-z_1^2 X_0/2
)\nonumber\\{}&&+f_{xx}^c(X_1+z_1X_0)^2/2+f_{xy}^c(X_1+z_1X_0)(Y_1+z_2Y_0)
+f_{yy}^c(Y_1+z_2Y_0)^2/2,
\end{eqnarray}
\begin{eqnarray}\label{Y2}
\dot Y_2(t)&=&A^sY_2(t)+z_2Y_1(t)+z_2^2f^s/2 +
f_x^s\big[X_2+z_1X_1+z_1^2X_0/2-z_2(X_1+z_1X_0)\big]+f_y^s(Y_2-z_2^2Y_0/2)\nonumber\\{}&&+
f_{xx}^s(X_1+z_1X_0)^2/2 +f_{xy}^s(X_1+z_1X_0)(Y_1+z_2Y_0)+
f_{yy}^s(Y_1+z_2Y_0)^2/2.
\end{eqnarray}
 Hence,
\begin{eqnarray*}
\tilde H^\varepsilon(\omega, \tilde\xi)&=&\int_{-\infty}^0
e^{-A^s\tau+\varepsilon\int_\tau^0
 z_2(\theta_r\omega)\,\mathrm{d}r}F^s\big(X(\tau), Y(\tau)\big)\,\mathrm{d}\tau\\ &=&
\int_{-\infty}^0 e^{-A^s\tau}\Big(1+\varepsilon\int_\tau^0
z_2(\theta_r\omega)\,\mathrm{d}r + \frac{\varepsilon^2(\int_\tau^0
z_2(\theta_r\omega)\,\mathrm{d}r)^2}{2!}+\cdots\Big)\Big\{f^s\big(X_0(\tau),
Y_0(\tau)\big)
\\{}&&+ \varepsilon \Big[f^s_x\big(X_0(\tau), Y_0(\tau)\big)\big(X_1 + z_1X_0\big) +
f^s_y\big(X_0(\tau), Y_0(\tau)\big)\big(Y_1 + z_2Y_0\big) -
z_2f^s\big(X_0(\tau),
Y_0(\tau)\big)\Big]\\{}&&+\varepsilon^2\Big[z_2^2f^s/2 +
f_x^s\big[X_2+z_1X_1+z_1^2X_0/2-z_2(X_1+z_1X_0)\big]+f_y^s(Y_2-z_2^2Y_0/2)\\{}&&+
f_{xx}^s(X_1+z_1X_0)^2/2 +f_{xy}^s(X_1+z_1X_0)(Y_1+z_2Y_0)+
f_{yy}^s(Y_1+z_2Y_0)^2/2\Big]\Big\}\,\mathrm{d}\tau
+\mathcal{O}(\varepsilon^3)\\&=& \int_{-\infty}^0
e^{-A^s\tau}f^s\big(X_0(\tau), Y_0(\tau)\big)\,\mathrm{d}\tau +
\varepsilon\Big\{\int_{-\infty}^0 e^{-A^s\tau}\Big[\int_\tau^0
z_2(\theta_r\omega)\,\mathrm{d}r \,f^s\big(X_0(\tau), Y_0(\tau)\big)
\\{}&& + f^s_x\big(X_0(\tau), Y_0(\tau)\big)\big(X_1 +
z_1X_0\big) + f^s_y\big(X_0(\tau), Y_0(\tau)\big)\big(Y_1 +
z_2Y_0\big) - z_2f^s\big(X_0(\tau),
Y_0(\tau)\big)\Big]\,\mathrm{d}\tau\Big\}\\{}&&+
\varepsilon^2\Big\{\int_{-\infty}^0 e^{-A^s\tau}\Big[z_2^2f^s/2 +
f_x^s\big[X_2+z_1X_1+z_1^2X_0/2-z_2(X_1+z_1X_0)\big]+f_y^s(Y_2-z_2^2Y_0/2)\nonumber\\{}&&+
f_{xx}^s(X_1+z_1X_0)^2/2 +f_{xy}^s(X_1+z_1X_0)(Y_1+z_2Y_0)+
f_{yy}^s(Y_1+z_2Y_0)^2/2+
(\int_\tau^0z_2\mathrm{d}r)^2f^s/2\nonumber\\{}&&
+\int_\tau^0z_2\mathrm{d}r
\Big(-z_2f^s+f_x^s(X_1+z_1X_0)+f_y^s(Y_1+z_2Y_0)\Big)
\Big]\,\mathrm{d}\tau\Big\}+\mathcal{O}(\varepsilon^3).
\end{eqnarray*}
Therefore, matching with \eqref{H=H^d+H^1}, we get
\begin{eqnarray}\label{t-H^d}
\tilde H^d(\tilde\xi)=\int_{-\infty}^0
e^{-A^s\tau}f^s\big(X_0(\tau), Y_0(\tau)\big)\,\mathrm{d}\tau,
\end{eqnarray}
\begin{eqnarray}\label{t-H^1}
\nonumber \tilde H^1(\omega, \tilde\xi)&=&\int_{-\infty}^0
e^{-A^s\tau}\Big\{\int_\tau^0 z_2(\theta_r\omega)\,\mathrm{d}r
\,f^s\big(X_0(\tau), Y_0(\tau)\big) + f^s_x\big(X_0(\tau),
Y_0(\tau)\big)\big[X_1(\tau) + z_1(\theta_\tau\omega)X_0(\tau)\big]
\\{}&& + f^s_y\big(X_0(\tau), Y_0(\tau)\big)\big[Y_1(\tau) +
z_2(\theta_\tau\omega)Y_0(\tau)\big] -
z_2(\theta_\tau\omega)f^s\big(X_0(\tau),
Y_0(\tau)\big)\Big\}\,\mathrm{d}\tau,
\end{eqnarray}
and
\begin{eqnarray}\label{t-H^2}
\nonumber \tilde H^2(\omega, \tilde\xi)&=&\int_{-\infty}^0
e^{-A^s\tau}\Big[z_2^2f^s/2 +
f_x^s\big[X_2+z_1X_1+z_1^2X_0/2-z_2(X_1+z_1X_0)\big]+f_y^s(Y_2-z_2^2Y_0/2)\nonumber\\{}&&+
f_{xx}^s(X_1+z_1X_0)^2/2 +f_{xy}^s(X_1+z_1X_0)(Y_1+z_2Y_0)+
f_{yy}^s(Y_1+z_2Y_0)^2/2+
(\int_\tau^0z_2\mathrm{d}r)^2f^s/2\nonumber\\{}&&
+\int_\tau^0z_2\mathrm{d}r
\Big(-z_2f^s+f_x^s(X_1+z_1X_0)+f_y^s(Y_1+z_2Y_0)\Big)
\Big]\,\mathrm{d}\tau.
\end{eqnarray}
Moreover, by \eqref{inversion}, the center manifold  for the
original stochastic system \eqref{stochastic center} is
$$\mathcal{H}^\varepsilon(\omega) =\big\{\big(\xi, H^\varepsilon(\omega, \xi
)\big) \big\}=\big\{\big(\xi, H^d(\xi)+\varepsilon H^1(\omega,
\xi)\varepsilon^2 H^2(\omega, \xi)++\mathcal{O}(\varepsilon^3)\big)
\big\},$$ with
\begin{eqnarray}\label{H^d}
 H^d(\xi)=\tilde H^d(\xi),
\end{eqnarray}
\begin{eqnarray}\label{H^1}
 H^1(\omega, \xi)=\tilde H^1(\omega, \xi)+
 z_2(\omega)\tilde H^d(\xi)-
 \tilde H_\xi^d(\xi)z_1(\omega)\xi,
\end{eqnarray}
and
\begin{eqnarray}\label{H^2}
 H^2(\omega, \xi)=\tilde H^d_\xi(\xi) z_1^2 \xi/2 -z_2 \tilde H^d_\xi(\xi) z_1
\xi + z_2^2\tilde H^d(\xi)/2 -\tilde H^1_\xi(\omega, \xi) z_1 \xi +
z_2\tilde H^1(\omega, \xi) + \tilde H^2(\omega, \xi).
\end{eqnarray}

\bigskip
We summarize the above approximation result in the following theorem.

\begin{theorem}[Approximation of a random center manifold] \label{center999}
 Assume that $\mathbf{(H)}$ holds and that $\eta$ satisfies $\gamma<\eta<\beta$ and $\frac{KL_f}{\eta-\gamma}+\frac{KL_g}{\beta-\eta}<1$. Then there
exists a   center manifold  $\mathcal{
H}^\varepsilon(\omega)=~\big\{\big(\xi,  H^\varepsilon(\omega, \xi
)\big)\big\}$ for the stochastic system \eqref{stochastic center},
where $H^\varepsilon(\omega, \xi)=H^d(\xi)+\varepsilon H^1(\omega,
\xi)+\varepsilon^2 H^2(\omega, \xi)+\mathcal{O}(\varepsilon^3)$,
with $H^d(\xi)$, $H^1(\omega, \xi)$ and $H^2(\omega, \xi)$   in
\eqref{H^d}, \eqref{H^1} and \eqref{H^2}, respectively.
\end{theorem}

\bigskip

We now look at an example of the approximation of random center manifolds.
\begin{example}
Consider a SDE system with Lyapunov exponent $\lambda_c=0$ and
$\lambda_s<0$
\begin{eqnarray}
\begin{cases}\label{eg_center}
\dot x = a^c x + \varepsilon x\circ \dot W_t,\,\,\, x\in\mathbb{R},\\
\dot y = -y -x^2 + \varepsilon y \circ \dot W_t,\,\,y\in\mathbb{R}.
\end{cases}
\end{eqnarray}
 Here $A^c=a^c$, $A^s=-1$, $f^c(x, y)=0$, $f^s(x,
y)=-x^2$, $z(\omega)=\int_{-\infty}^0e^\tau\,\mathrm{d}W_\tau$ and
$z(\theta_t\omega)=\int_{-\infty}^t
e^{\tau-t}\,\mathrm{d}W_\tau=e^{-t}z(\omega)+e^{-t}\int_{0}^te^\tau\,\mathrm{d}W_\tau$.
We restrict this system on a bounded disk containing the origin. The
new system satisfies the assumption \textbf{(H)}. In this way, we
obtain a local center manifold and we now consider its small noise
approximation.

 We transform  this SDE system into a RDE system,
\begin{eqnarray}
\begin{cases}
\dot X(t) = a^c X + \varepsilon z(\theta_t\omega) X,\,\,\, x\in\mathbb{R},\\
\dot Y(t) = -Y  + \varepsilon z(\theta_t\omega) Y -e^{-\varepsilon
z(\theta_t\omega)}(e^{\varepsilon z(\theta_t\omega)}X)^2
,\,\,y\in\mathbb{R}.
\end{cases}
\end{eqnarray}

Denote $a=a^c + \varepsilon z(\theta_t\omega)$, $b=-1+\varepsilon
z(\theta_t\omega)$, and $A=\begin{pmatrix}  a \quad 0 \\ 0\quad b
\\\end{pmatrix}$.  By   the property\\
$\lim\limits_{t\rightarrow
\pm\infty}\frac{1}{t}\int_0^tz(\theta_\tau
\omega)\,\mathrm{d}\tau=0$ in \cite{Duan},  we have the Lyapunov
exponents
\begin{eqnarray}
\begin{cases}\lambda_c=\lim\limits_{t\rightarrow
\pm\infty}\,\,\frac{1}{t}\log | e^{\int_0^t
A(\tau)\,\mathrm{d}\tau}u |=a^c=0,\quad \text{when}\,\,\,
u=\begin{pmatrix}  u_1
\\ 0\\\end{pmatrix},\\ \lambda_s=\lim\limits_{t\rightarrow
\pm\infty}\,\,\frac{1}{t}\log | e^{\int_0^t
A(\tau)\,\mathrm{d}\tau}u |=-1<0, \quad\text{when}\,\,\,
u=\begin{pmatrix}  0
\\ u_2\\\end{pmatrix}.
\end{cases}
\end{eqnarray}
By \eqref{X0Y0} and \eqref{t-H^d}, we get\\
$$X_0(t)=\tilde\xi,\quad\quad\tilde H^d(\tilde\xi)=\int_{-\infty}^0
e^\tau(-X_0^2)\,\mathrm{d}\tau=-\tilde\xi^2,$$ and
$$Y_0(t)=e^{-t}\tilde H^d(\tilde\xi) -
\tilde\xi^2(1-e^{-t})=-\tilde\xi^2.$$ From \eqref{X1Y1}, \eqref{X2},
\eqref{t-H^1} and \eqref{t-H^2}, we obtain
$$X_1(t)=X_1(0)+\tilde\xi\int_0^tz(\theta_\tau\omega)\,\mathrm{d}\tau=\tilde\xi\int_0^tz(\theta_\tau\omega)\,\mathrm{d}\tau,$$
$$X_2(t)=X_2(0)+\tilde\xi\int_0^tz(\theta_s\omega)\int_0^sz(\theta_r\omega)\,\mathrm{d}r\,\mathrm{d}s
=\tilde\xi\int_0^tz(\theta_s\omega)\int_0^sz(\theta_r\omega)\,\mathrm{d}r\,\mathrm{d}s,$$
and\begin{eqnarray} \nonumber\tilde H^1(\omega,
\tilde\xi)&=&\int_{-\infty}^0e^\tau\big[\int_\tau^0z(\theta_r\omega)\,dr\big(-X_0^2(\tau)\big)-2X_0(\tau)\big(X_1(\tau)
+z(\theta_\tau\omega)X_0(\tau)\big)+z(\theta_\tau\omega)X^2_0(\tau)\big]\,\mathrm{d}\tau\\
&=&\nonumber
-\tilde\xi^2\int_{-\infty}^0e^\tau\big[\int_\tau^0z(\theta_r\omega)\,\mathrm{d}r
+ 2\big(\int_0^\tau z(\theta_r\omega)\,\mathrm{d}r +
z(\theta_\tau\omega)\big) -
z(\theta_\tau\omega)\big]\,\mathrm{d}\tau,
\end{eqnarray}
\begin{eqnarray} \nonumber\tilde H^2(\omega,
\tilde\xi)&=&\int_{-\infty}^0e^\tau\big[-z^2\xi^2/2-2\xi(X_2-z^2X_0/2)-(X_1+zX_0)^2
-\xi^2(\int_\tau^0z(\theta_r\omega)\,\mathrm{d}r)^2/2\\{}&&\quad\quad\quad\quad+\int_\tau^0z(\theta_r\omega)\,\mathrm{d}r
\big(z\xi^2-2\xi(X_1+zX_0)\big)\big]\,\mathrm{d}\tau\nonumber\\
&=&\nonumber
\tilde\xi^2\int_{-\infty}^0e^\tau\big[-\frac{z^2}{2}-2\int^\tau_0z\int_0^s
z(\theta_r\omega)\,\mathrm{d}r \,\mathrm{d}s+ \frac{\big(\int_0^\tau
z(\theta_r\omega)\,\mathrm{d}r\big)^2}{2} +
z(\theta_\tau\omega)\int_\tau^0z\,\mathrm{d}r\big]\,\mathrm{d}\tau,
\end{eqnarray}
by calculation
\begin{eqnarray}
\tilde H^1(\omega,
\tilde\xi)&=&\tilde\xi^2\big[\int_{-\infty}^0e^\tau\int_\tau^0z(\theta_r\omega)\,\mathrm{d}r\,\mathrm{d}\tau
-\int_{-\infty}^0e^\tau z(\theta_\tau\omega)\,\mathrm{d}\tau\big]=0,
\end{eqnarray}
\begin{eqnarray}
\tilde H^2(\omega,
\tilde\xi)&=&-4\tilde\xi^2\int_{-\infty}^0\int_{-\infty}^ve^u\,\mathrm{d}W_u\,\mathrm{d}W_v+
\tilde\xi^2(\int_{-\infty}^0e^u\,\mathrm{d}W_u)^2
+2\tilde\xi^2\int_{-\infty}^0\int_v^0e^v(u-v)\,\mathrm{d}W_u\,\mathrm{d}W_v.\nonumber\\
\end{eqnarray}
By \eqref{H^d}, \eqref{H^1} and \eqref{H^2} $$H^d(\xi)=\tilde
H^d(\xi)=-\xi^2,$$
\begin{eqnarray*}
H^1(\omega, \xi)&=&\tilde H^1(\omega, \xi) + z(\omega)\tilde
H^d(\xi) - \tilde H_\xi^d(\xi)z(\omega)\xi
=\xi^2\int_{-\infty}^0e^\tau\,\mathrm{d}W_\tau,
\end{eqnarray*}and
\begin{eqnarray}
 H^2(\omega, \xi)&=&\tilde H^d_\xi(\xi) z_1^2 \xi/2 -z_2 \tilde H^d_\xi(\xi) z_1
\xi + z_2^2\tilde H^d(\xi)/2 -\tilde H^1_\xi(\omega, \xi) z_1 \xi +
z_2\tilde H^1(\omega, \xi) + \tilde H^2(\omega,
\xi)\nonumber\\&=&\frac{z^2\xi^2}{2}-4\xi^2\int_{-\infty}^0\int_{-\infty}^ve^u\,\mathrm{d}W_u\,\mathrm{d}W_v+
\xi^2(\int_{-\infty}^0e^u\,\mathrm{d}W_u)^2
+2\xi^2\int_{-\infty}^0\int_v^0e^v(u-v)\,\mathrm{d}W_u\,\mathrm{d}W_v.\nonumber
\end{eqnarray}

Therefore, we have an approximation for   center manifold of
\eqref{eg_center} for $\varepsilon$ sufficiently small

$$\mathcal{H}^\varepsilon(\omega) =\big\{\big( \xi,
H^\varepsilon(\omega, \xi )\big)\big\},$$
with
\begin{eqnarray*}
H^\varepsilon(\omega, \xi
)&=&-\xi^2+\varepsilon\,\xi^2\int_{-\infty}^0e^\tau\,dW_\tau+
\varepsilon^2
\xi^2\Big[\frac{z^2}{2}-4\int_{-\infty}^0\int_{-\infty}^ve^u\,\mathrm{d}W_u\,\mathrm{d}W_v\\{}&&+
(\int_{-\infty}^0e^u\,\mathrm{d}W_u)^2
+2\int_{-\infty}^0\int_v^0e^v(u-v)\,\mathrm{d}W_u\,\mathrm{d}W_v\Big]+\mathcal{O}(\varepsilon^3).
\end{eqnarray*}
The center manifold of the corresponding deterministic   system
\begin{eqnarray}
\begin{cases}
\dot x = 0,\\
\dot y = -y - x^2 .
\end{cases}
\end{eqnarray}
is the graph of the function $h^c(\xi)=-\xi^2$. Realizations of the
approximate center manifold with error order 2 or 3 in
$\varepsilon$, together with the deterministic center manifold, are
shown in Figure \ref{center02}.
\begin{figure}
\includegraphics[height=7cm,width=15cm]{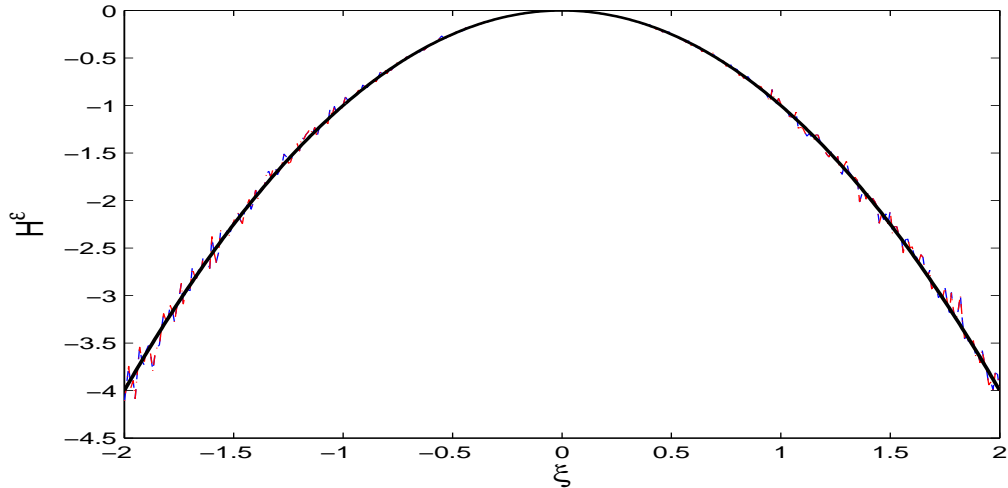}
\caption{Approximate center manifold  of system \eqref{eg_center}:
 $\varepsilon=0$ (no noise,  black or solid curve) v.s.
$\varepsilon=0.05$ (Red or dash-dot curve for
$H^\varepsilon=H^d+\varepsilon H^1$, blue or dash curve for
$H^\varepsilon=H^d+ \varepsilon H^1+\varepsilon^2H^2$). }
\label{center02}
\end{figure}
\end{example}


\medskip

\textbf{Acknowledgements}. We would like to thank Xiaopeng Chen and Jinqiao Duan for helpful discussions.

\end{document}